\documentclass{article}
\usepackage{amsfonts}
\usepackage{amsmath}
\usepackage{xypic}
\usepackage{graphicx,epsfig}
\usepackage{amssymb}
\usepackage{rotating}
\usepackage{theorem}
\newtheorem{theo}{Theorem}[section]
\newtheorem{cor}{Corollary}[section]
\theorembodyfont{\rmfamily}
\usepackage{float}
\usepackage[bf]{caption}

\begin{document}
\thispagestyle{plain}
\begin{center}
\begin{Large}
\textbf{Universal Cone Manifolds and the Poincar\'e Conjecture I}\\
\vspace{2cm}
Ali Aalam
\end{Large}
\end{center}
\vspace{1cm}
\textbf{Abstract} We identify a universal group $U$ and show that $\Bbb H^3/G$ is $S^3$ when $G$ is a finite index subgroup of $U$ generated by elements of finite order.

\section*{}In \cite{Aal} we described a family of hyperbolic cone manifolds $W^3$ whose singular locus is the Whitehead link. We now consider the case where the isotropy groups are cyclic of orders 4 and 8 on suitable link components. Thus, this orbifold is defined by a finitely generated group $U$ of orientation-preserving isometries of hyperbolic three-space $\Bbb H^3$.

We show that $U$ is {\it universal}, i.e., for every closed and oriented three-manifold $M^3$ there exists a subgroup $G\leq U$ of finite index such that $M^3$ is homeomorphic to the quotient $\Bbb H^3/G$. Thus $M^3$ is a hyperbolic orbifold covering $W^3$. A {\it universal cone manifold} is a cone manifold whose singular set is a universal knot or link.

Every subgroup $G \leq U$ of finite index of a universal group $U$ defines a closed, oriented three-manifold $\Bbb H^3/G$ whose fundamental group is $G/F$, where $F$ is the subgroup of $G$ generated by the elements of finite order. Therefore, the simply connected three-manifolds correspond precisely to the subgroups $G\leq U$ of finite index that are generated by elements of finite order. This offers an approach to the Poincar\'e Conjecture.

Following \cite{HLM4} we show that every closed, oriented three-manifold $M^3$ is a covering of $S^3$ branched along the Whitehead link with branching indices $1,2,4$ and $8$. Universality of the Whitehead link was proved in \cite{HLM1}, here we sharpen it by giving the branching indices.

We construct the {\it associated regular covering} branched over the Whitehead link and show that its branching indices are $4$ and $8$ on appropriate link components. We describe our construction of a hyperbolic orbifold as in \cite{Aal} which has the Whitehead link as singular locus and cyclic isotropy groups of orders $4$ and $8$. 
\section{Branching over the Whitehead link}
If $f \colon M^3\rightarrow S^3$ is a $d$-fold covering branched over a link $L\subset S^3$, the pre-image of a meridian $m$ consists of components which map onto $m$ as unbranched coverings of degrees $i_1,\dots ,i_k$. These numbers are the {\it brancing indices} of the component of $L$ with meridian $m$. 

In \cite{HLM1} it was shown that every closed, oriented 3-manifold is a covering of $S^3$ branched along the Whitehead link. We sharpen this results as follows.

\begin{theo}
Every closed, oriented 3-manifold is a covering of $S^3$ branched along the Whitehead link such that the branching indices of the components are (1,2,4) and (4,8).
\end{theo}

\begin{figure}[H]
\includegraphics[height=25mm,width=65mm]{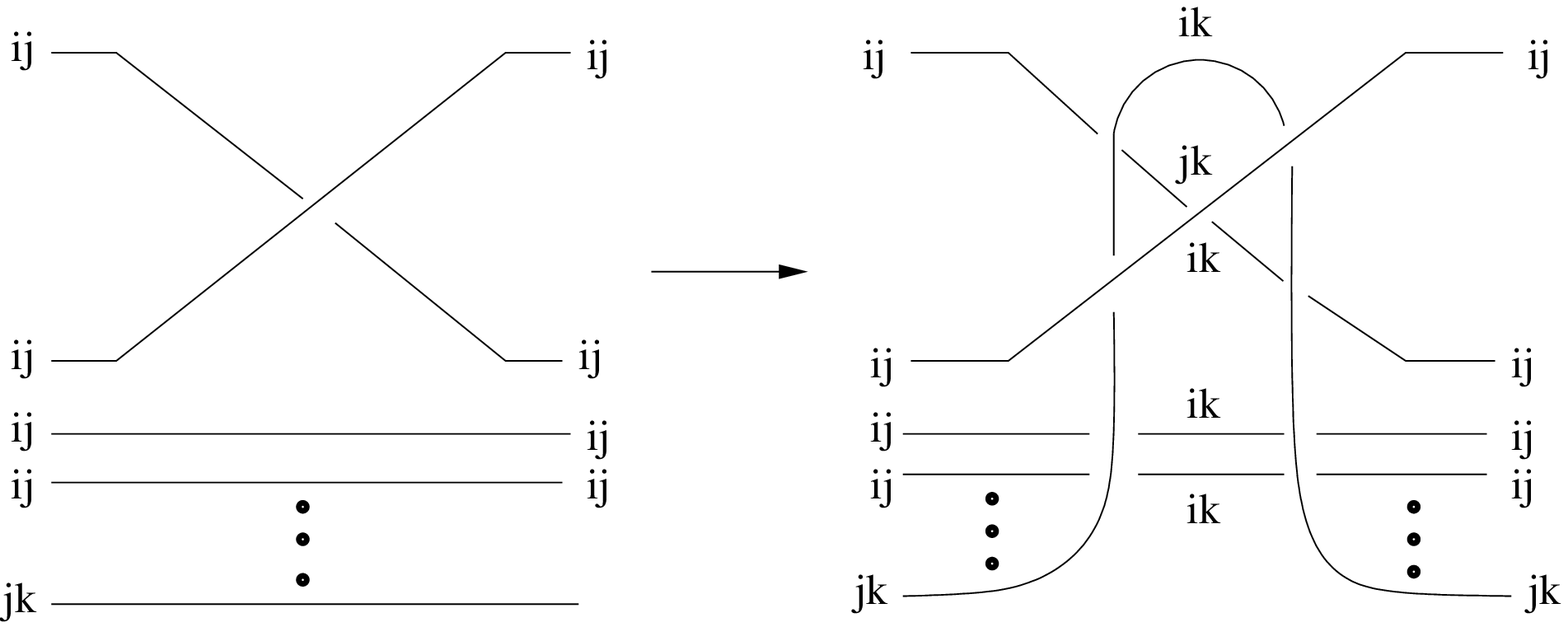}
\centering
\caption{}\label{fig:fig1}
\end{figure}

{\it Proof}. Let $M^3$ be a closed and oriented 3-manifold. Then there exists a 3-fold simple covering $f \colon M^3\rightarrow S^3$ branched over a knot $K$ (\cite{Hi}, \cite{Mon1}). The covering $f$ corresponds to a transitive representation

\begin{equation}\notag
\omega \colon \pi_1 (S^3 - K) \rightarrow S_3
\end{equation}
\\
onto the symmetric group of order 3, sending meridians to transpositions. We will prescribe $\omega$ by assigning transpositions $(1,2),(1,3),(2,3)$ to the overpasses of a normal projection of $K$. In each crossing the three intervening transpositions must be equal or pairwise different. We can further assume that $K$ is a closed braid, and that the transpositions intervening in a crossing are pairwise disjoint: it is sufficient to use the move shown in Figure \ref{fig:fig1}.

\begin{figure}[H]
\includegraphics[height=18mm,width=85mm]{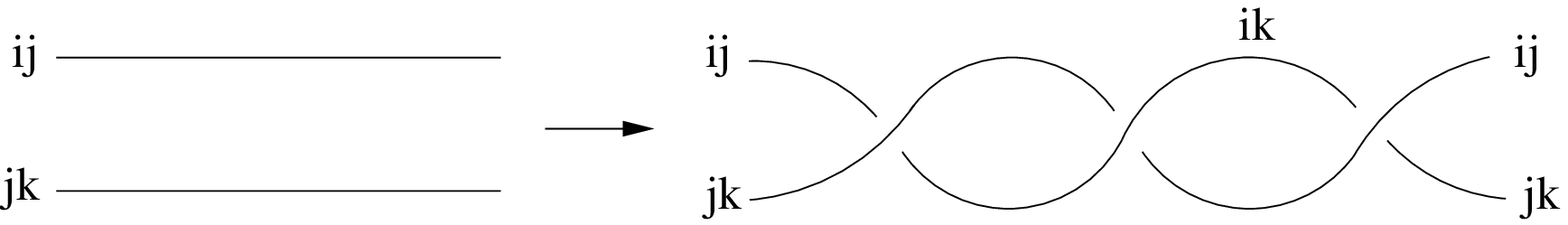}
\centering
\caption{}\label{fig:fig2}
\end{figure}

By using the move of Figure \ref{fig:fig2},

which does not alter $M^3$ \cite{Mon2} (see \cite{Mon3} and \cite{Fox}), \cite{Hir}, we can convert positive crossings into negative ones. Thus we have that $K$ is a negative closed braid such that the three transpositions intervening in each crossing are pairwise different.

The next operation is illustrated in Figure \ref{fig:fig3}. Here each crossing is converted into the boundary of an annulus together with a circle linking it. After applying this operation :

\begin{figure}[H]
\includegraphics[height=55mm, width=75mm]{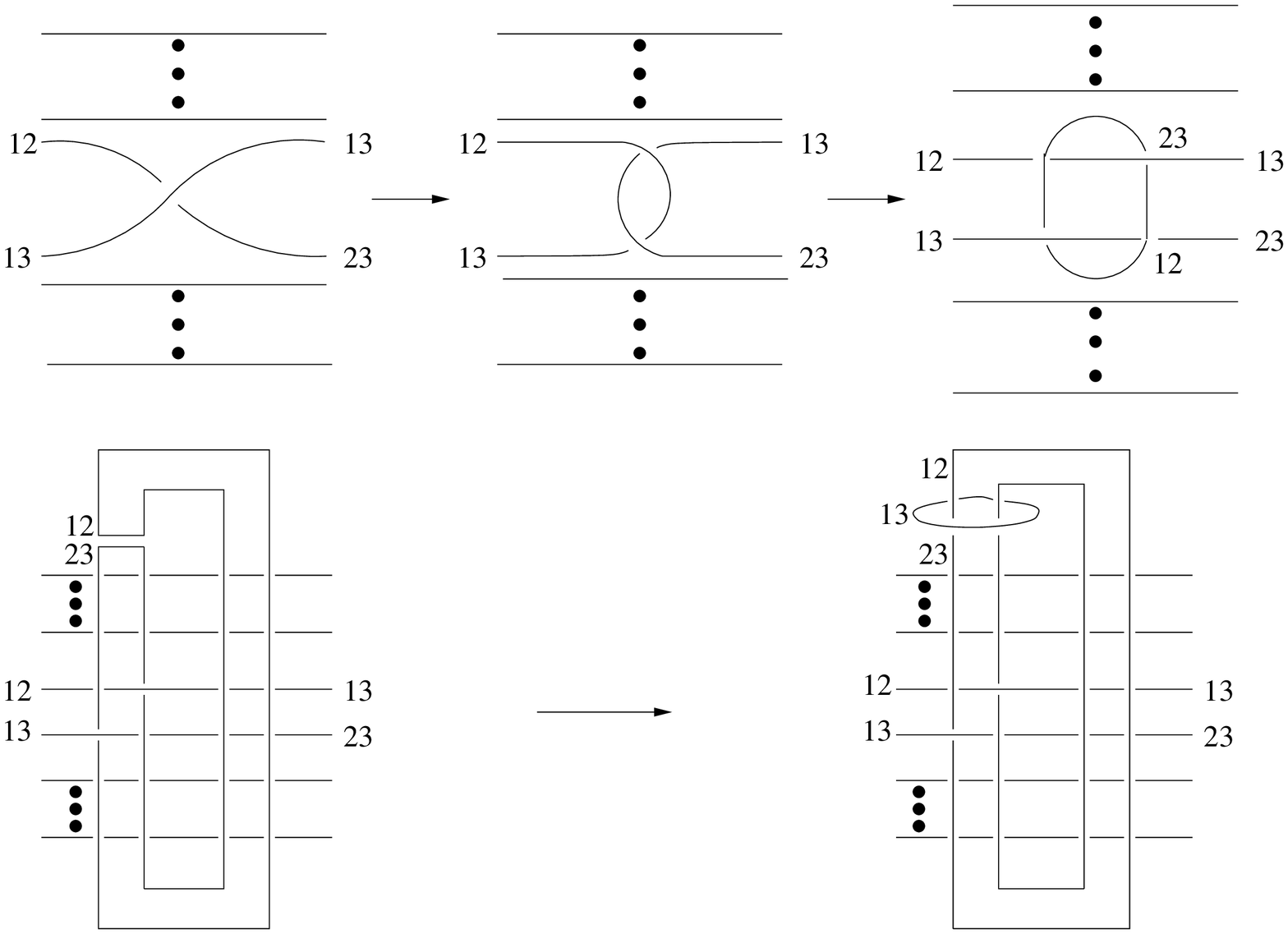}
\centering
\caption{}\label{fig:fig3}
\end{figure}

to all the crossings of $K$ we will have a number of annuli and circles linking them. These circles can be deformed into new horizontal strings as shown in Figure \ref{fig:fig4}.
\begin{figure}[H]
\includegraphics[height=32mm, width=75mm]{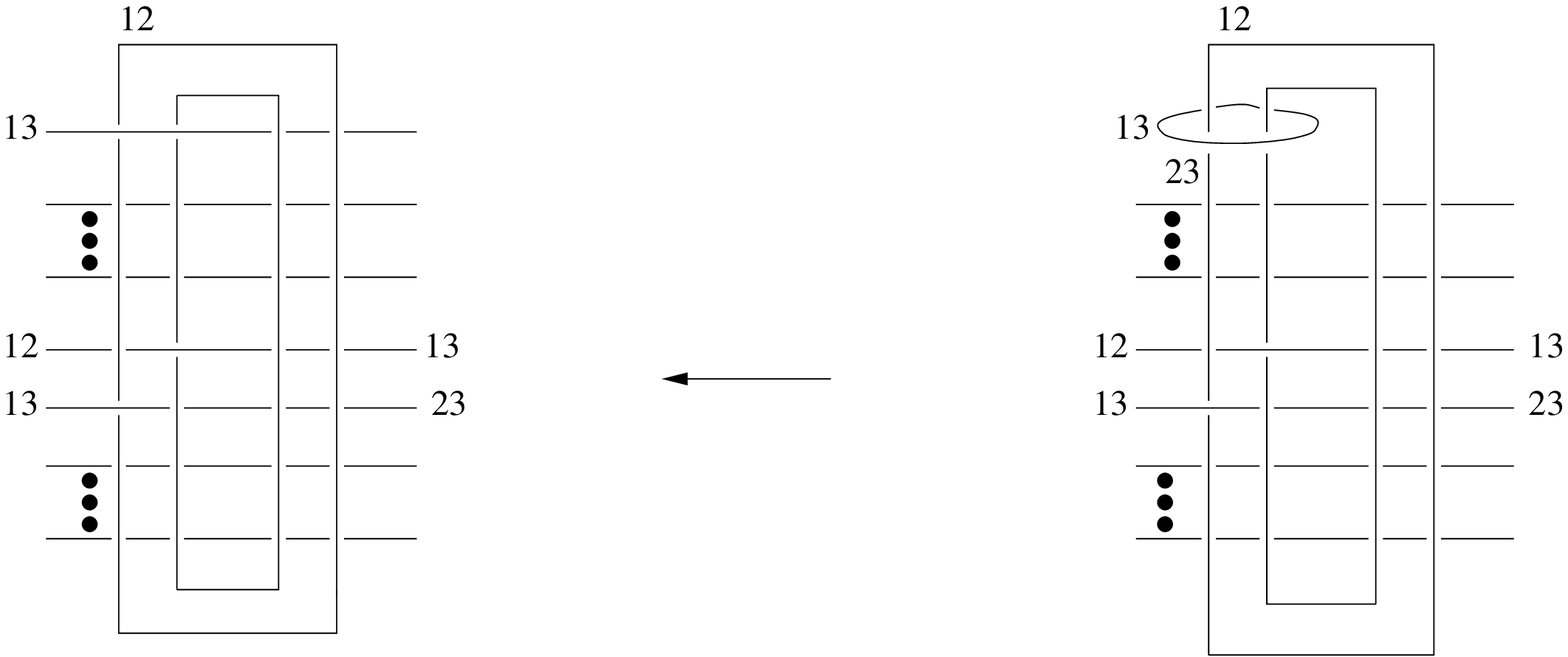}
\centering
\caption{}\label{fig:fig4}
\end{figure}

\begin{figure}[H]
\includegraphics[height=50mm, width=35mm]{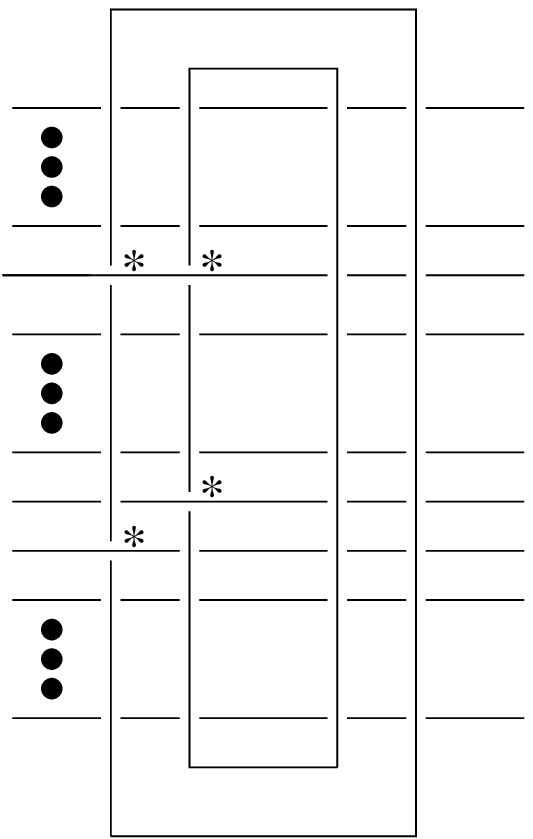}
\centering
\caption{}\label{fig:fig5}
\end{figure}

Therefore, we can assume that $K$ is composed by horizontal strings and annuli of the form shown in Figure \ref{fig:fig5}.

and since the covering is still a 3-fold simple covering, the branching indices of the components of $K$ are $(1,2)$.

For the next step of the proof, we need to recall the concept of {\it adding a trivial sheet}. In general, if $f \colon M^3 \rightarrow S^3$ is a $d$-fold branched covering, we obtain a $(d+1)$-fold branched covering $f' \colon M^3 \rightarrow S^3$ by enlarging the branching set $K$ of $f$ with a trivial knot. This knot must bound a disk not meeting $K$, and must be endowed with a transposition $(i,d+1), i<d+1$.
In general, if this trivial knot is endowed with the permutation $(i_1,d+1)(i_2,d+2) \dots (i_s,d+s),$ where $i_1,i_2,\dots ,i_s$ are all different and less than $d$, the result is a new branched covering $f':M^3 \rightarrow S^3$ with $(d+s)$ sheets.

\begin{figure}[H]
\includegraphics[height=25mm, width=115mm]{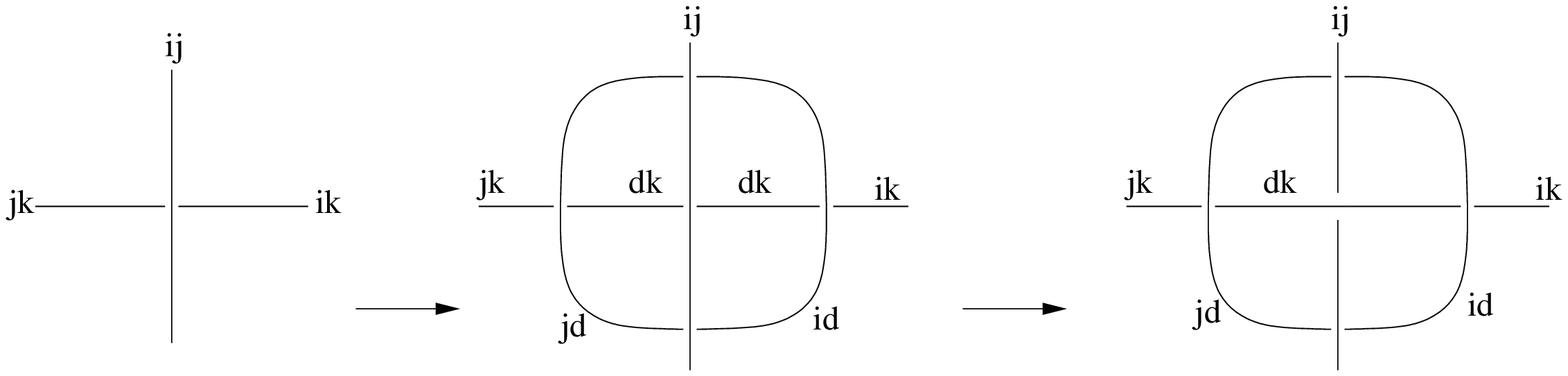}
\centering
\caption{}\label{fig:fig6}
\end{figure}

We now want the two left vertical strings of each annulus (Figure \ref{fig:fig5}) to go underneath the set of horizontal strings. To acheive this we apply two new moves to the crossings lying in the two left vertical strings of Figure \ref{fig:fig5}, and not endowed with stars. The first move converts an overcrossing to an undercrossing when the transpositions involved in the crossing are all different, this move is shown in Figure \ref{fig:fig6}

\begin{figure}[H]
\includegraphics[height=25mm, width=80mm]{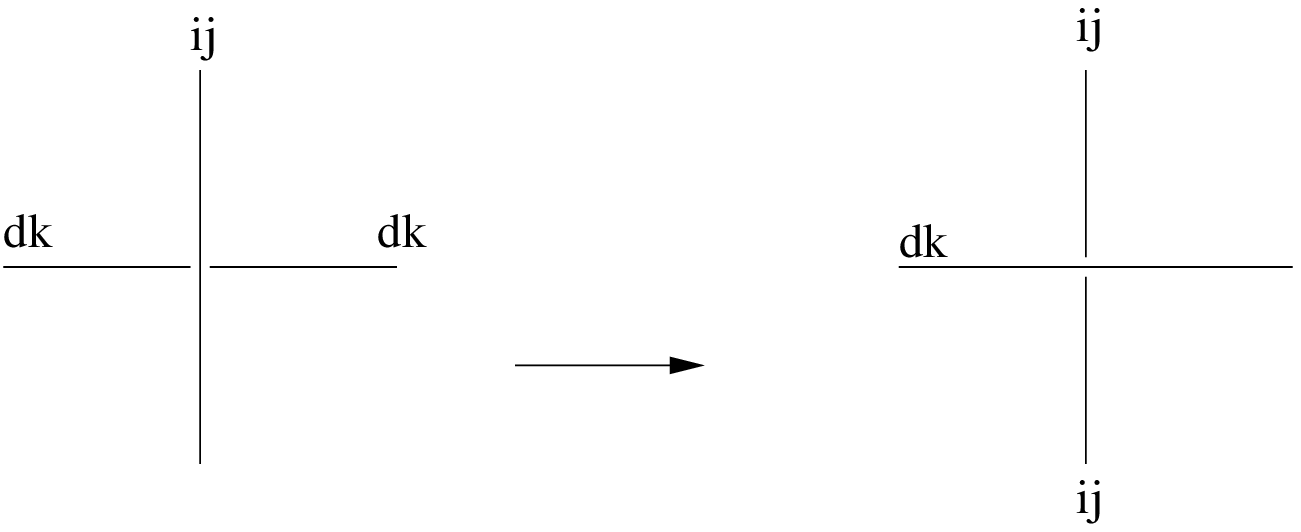}
\centering
\caption{}\label{fig:fig7}
\end{figure}

and it is the consequence of adding a new sheet and the move of Figure \ref{fig:fig7}

\begin{figure}[H]
\includegraphics[height=25mm, width=100mm]{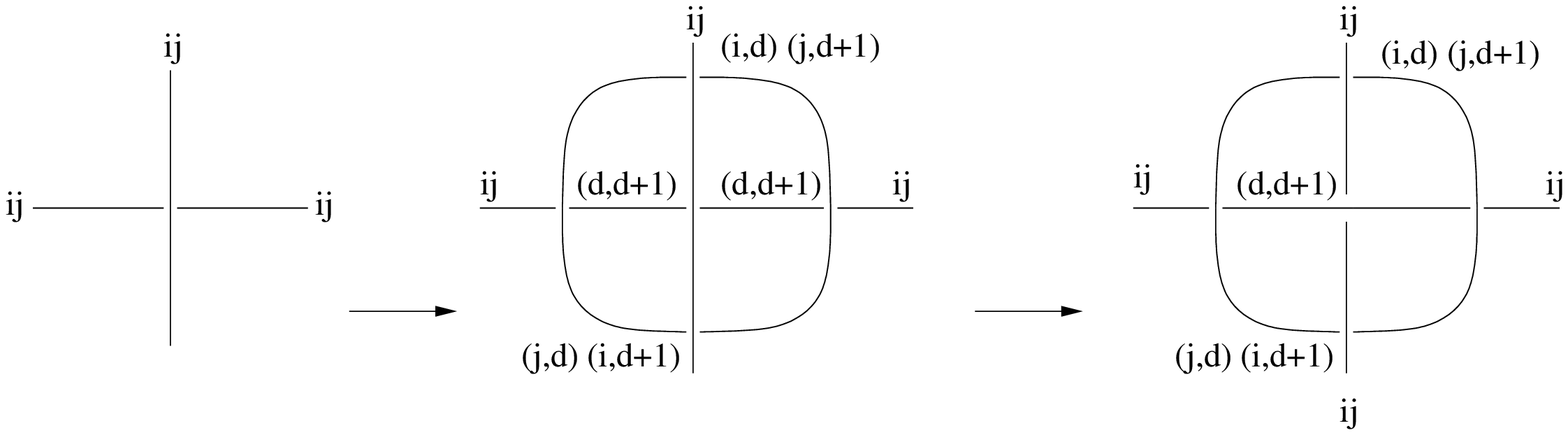}
\centering
\caption{}\label{fig:fig8}
\end{figure}

that does not change the covering manifold \cite{Mon2} and \cite{Mon3}. The second move produces the same effect as the first, i.e. changes an overcrossing into an undercrossing, but this time the transpositions involved in the crossing are all equal; this move is depicted in Figure \ref{fig:fig8}, and is similar to the move of Figure \ref{fig:fig6}; the difference consists of the addition of two new sheets to the covering.

After performing these moves in the relevant crossings, we have a branched covering $f_1:M^3 \rightarrow S^3$ which is not simple in general, but the branching indices for each component of the new branching set are $(1,2)$.

\begin{figure}[H]
\includegraphics[height=45mm, width= 25mm]{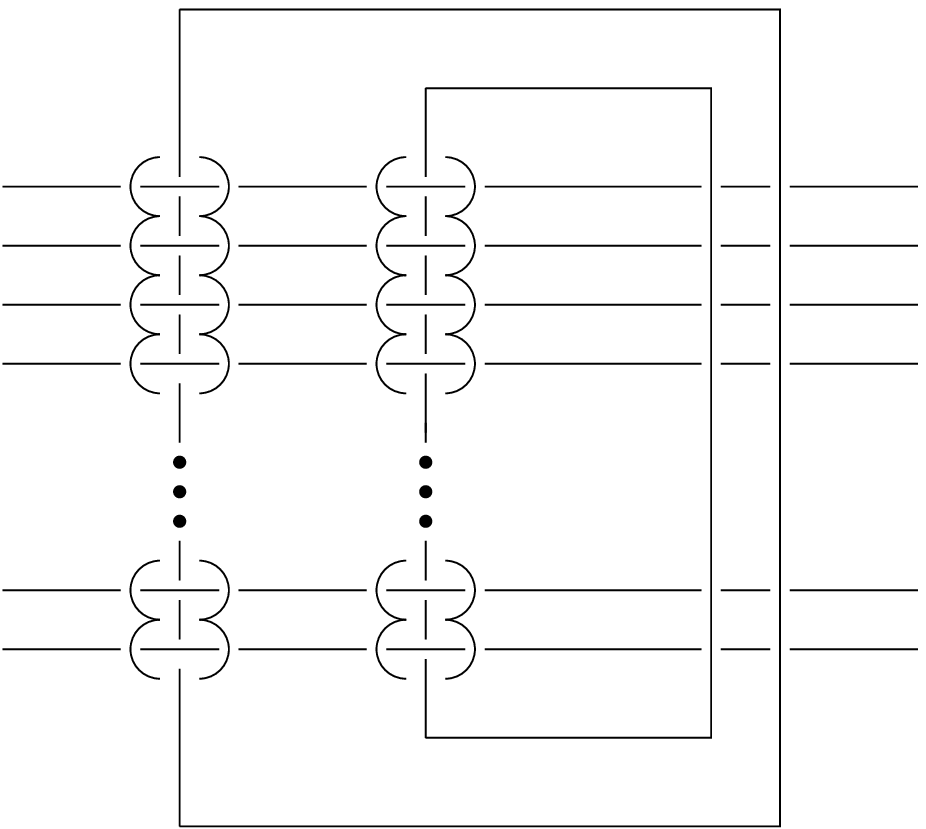}
\centering
\caption{}\label{fig:fig9}
\end{figure}

We can enlarge this new branching set by adding new curves with branching indices $(1)$, in such a way that each annulus has the form shown in Figure \ref{fig:fig9}.

\begin{figure}[H]
\includegraphics[height=65mm, width=65mm]{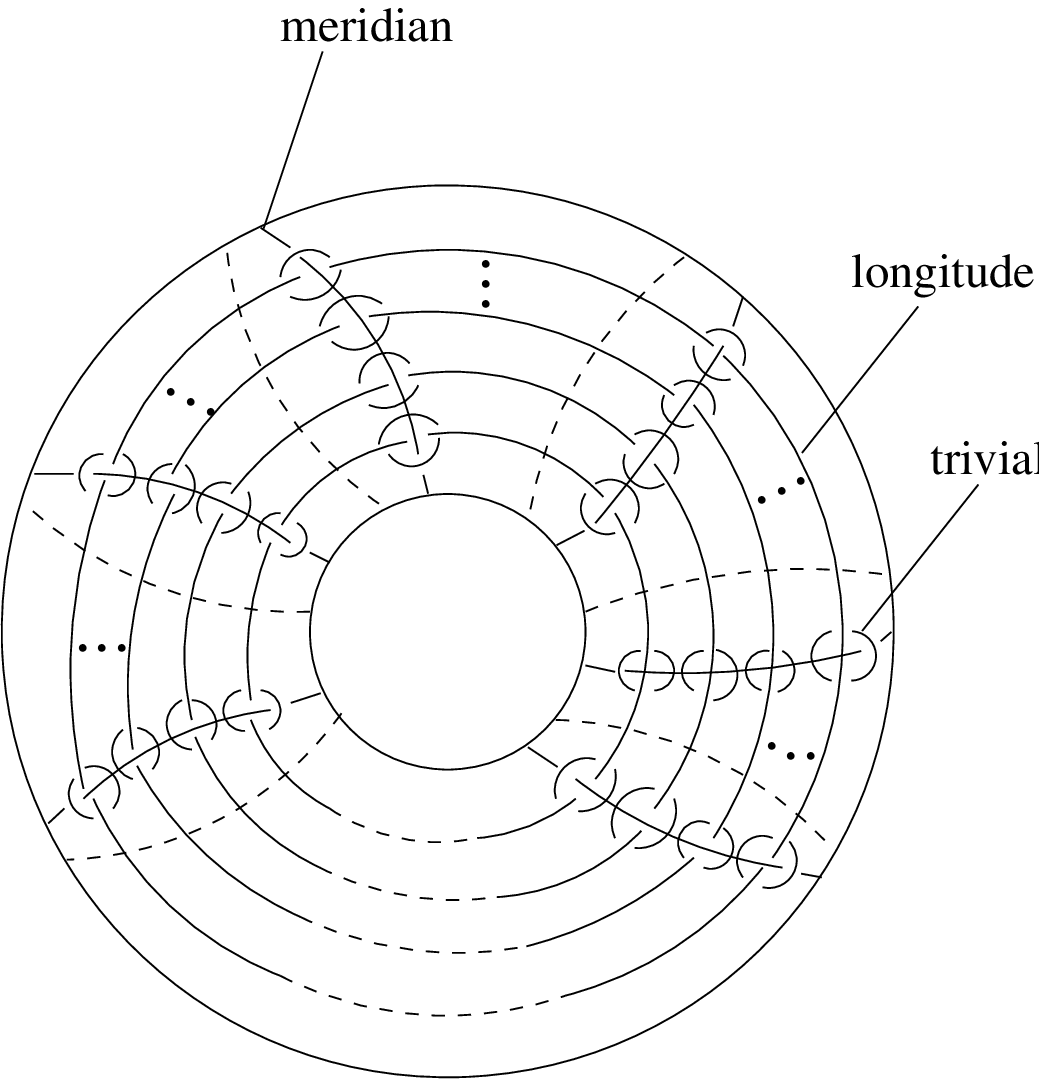}
\centering
\caption{}\label{fig:fig10}
\end{figure}

Thus we can assume the branching set $K_1$ of $f_1$ has a normal projection in a torus $T^2$ as shown in Figure \ref{fig:fig10}. Here the ``meridian'' and ``longitude'' components have branching indices $(1,2)$ and some others (the ones added last) have index $(1)$.

Furthermore, we can place the components of $K_1$ in such a way that the action of $\Bbb Z_m \times \Bbb Z_l$ in $S^3$ generated by a $2\pi/m$-rotation $M$ around the core $u$ of the unbounded solid torus, and by a $2\pi/l$-rotation $L$ around the core $b$ of the bounded solid torus, permutes the $m$ meridians, the $l$ longitudes and the trivial components (this idea of ``symmetrising with respect to a torus'' is due to Thurston \cite{Th1}).

\begin{figure}[H]
\includegraphics[height=18mm, width=64mm]{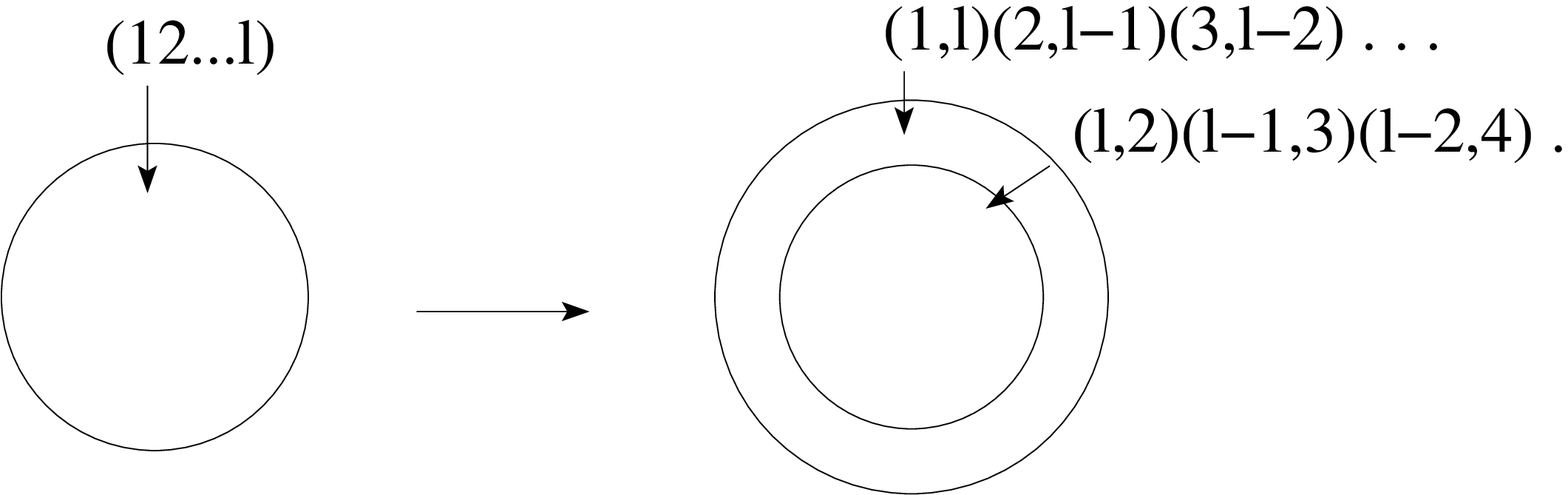}
\centering
\caption{}\label{fig:fig11}
\end{figure}

The $2\pi/l$-rotation $L$ generates the group of covering transformations of an $l$-fold cyclic covering $f_L \colon S^3 \rightarrow S^3$, whose branching set is $f_L(b)$. We can modify $f_L$ into $f'_L \colon S^3 \rightarrow S^3$, branched along two curves parallel to $f_L(b)$ and which do not link each other. The modification takes place in a regular neighbourhood $N^3$ of $f_L(b)$ outside of which the covering remains cyclic. To see this note that the monodromy $\eta_L \colon \pi_1(S^3 -f_L(b)) \rightarrow S_l$ sends a meridian of the branching set $f_L(b)$ of $f_L$ to the cycle $(12 \dots l)$. We modify this branching set and the monodromy as indicated in Figure \ref{fig:fig11} (see \cite{Mon3} for more details).

\begin{figure}[H]
\includegraphics[height=65mm, width=65mm]{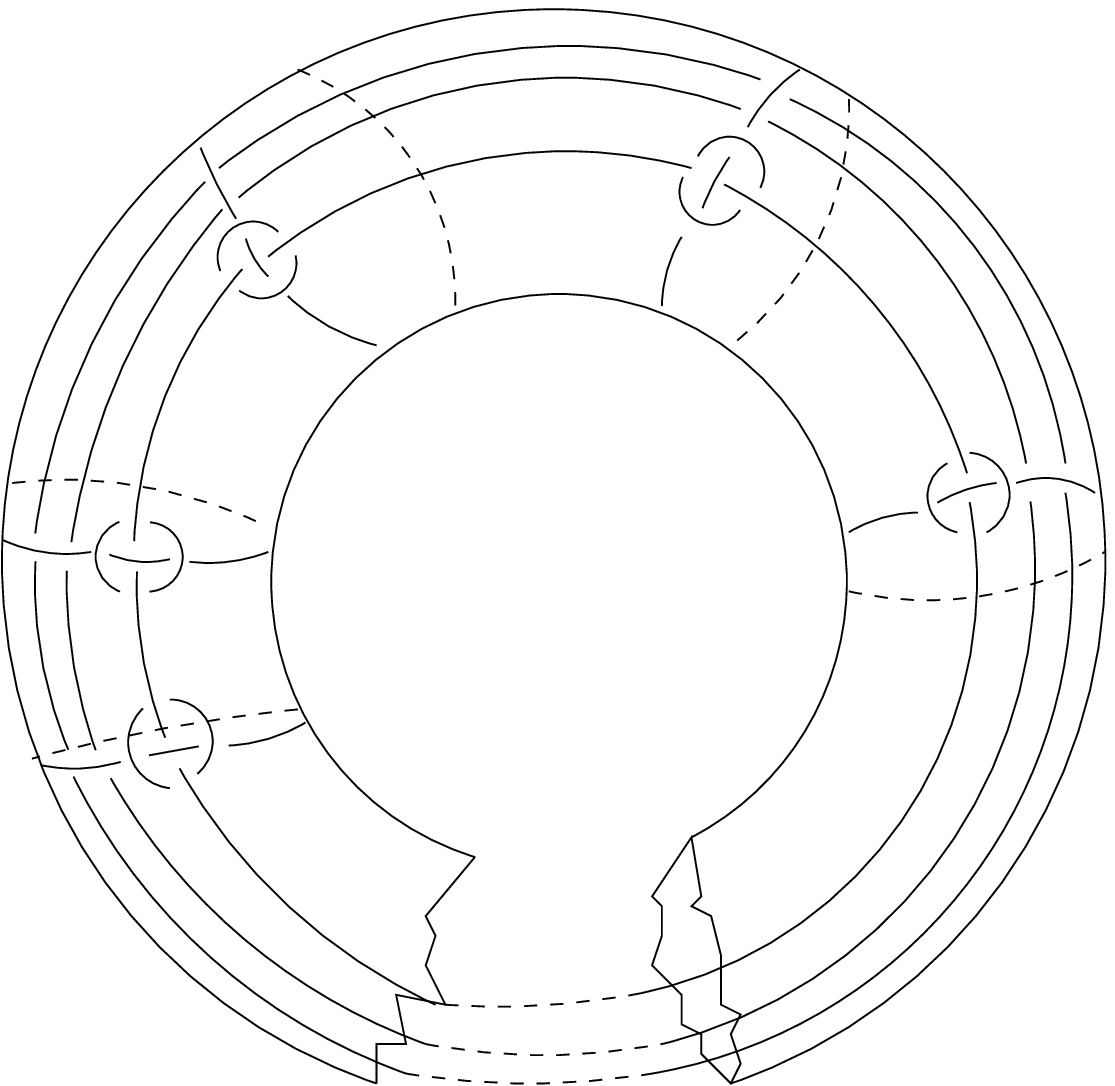}
\centering
\caption{}\label{fig:fig12}
\end{figure}

The image of $K_1$ under $f'_L$ together with the branching set of $f'_L$ form the branching set $K_2$ of the composition $f'_Lf_1 :=f_2 : M^3 \rightarrow S^3$. The link $K_2$ is depicted in Figure \ref{fig:fig12}. Each component has branching indices $(1,2)$.

\begin{figure}[H]
\includegraphics[height=35mm, width=115mm]{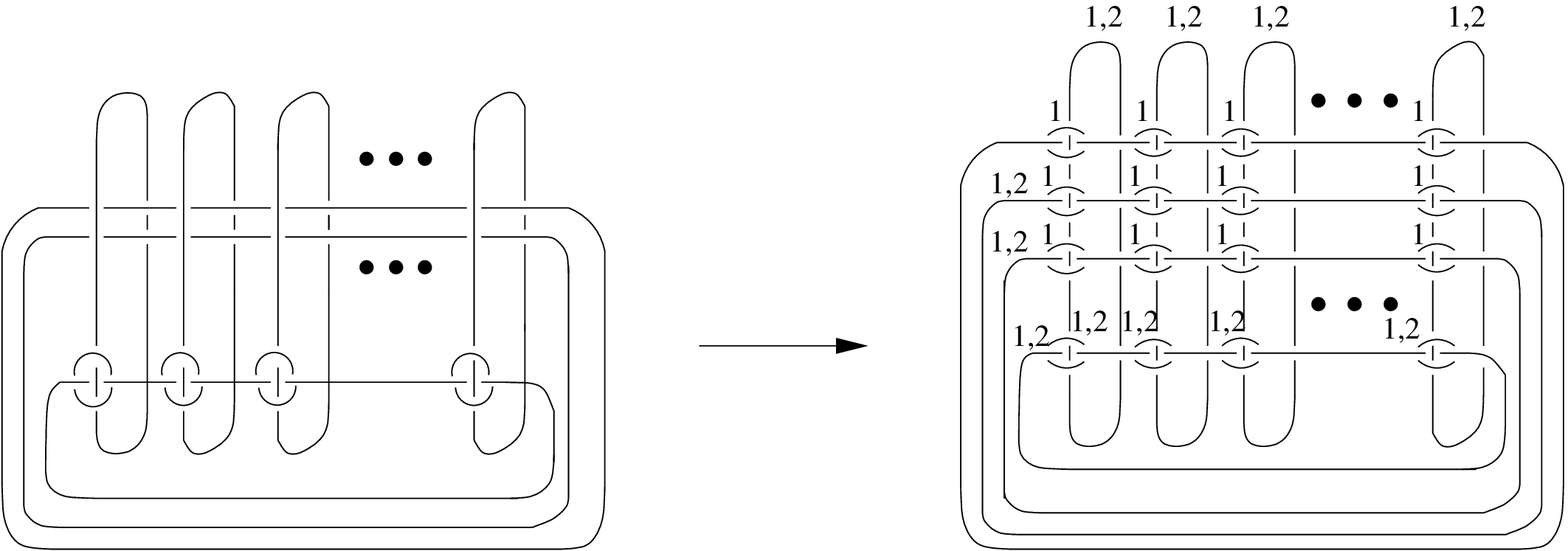}
\centering
\caption{}\label{fig:fig13}
\end{figure}

We now enlarge $K_2$ by adding a new curve, parallel to the two in the bounded solid torus, and $3m$ trivial curves. The new curves have branching index $(1)$ (this process is schematically detailed in Figure \ref{fig:fig13}). The result is the covering $f_2 : M^3 \rightarrow S^3$ branched over $K'_2$, as in Figure \ref{fig:fig10} but with $l=4$.

\begin{figure}[H]
\includegraphics[height=30mm, width=55mm]{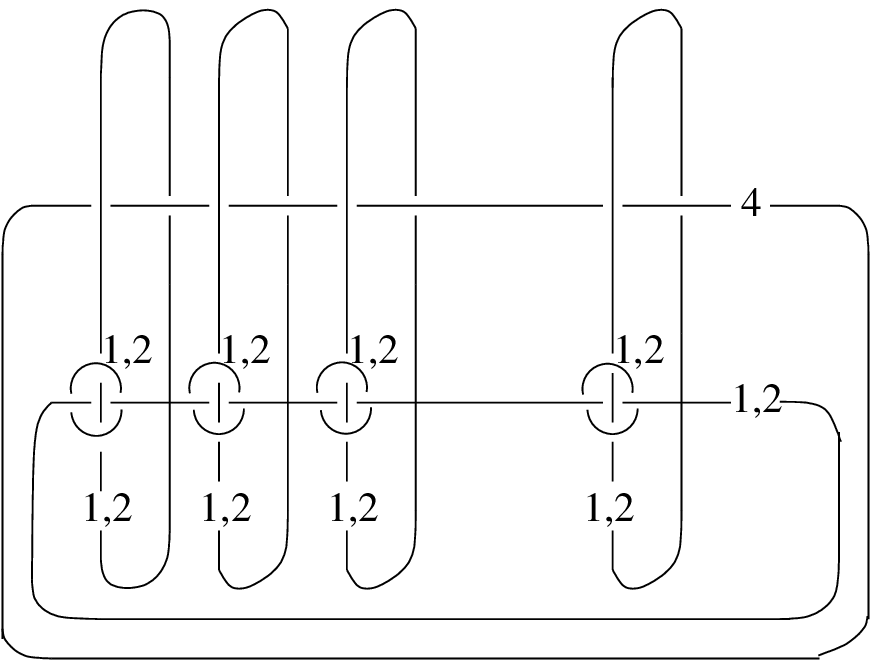}
\centering
\caption{}\label{fig:fig14}
\end{figure}

We compose $f_2$ with $f_L \colon S^3 \rightarrow S^3$ for $l=4$ and we obtain $f_3 : M^3 \rightarrow S^3$ with the branching set $K_3$ shown in Figure \ref{fig:fig14}. Here, all components have branching indices $(1,2)$ except the core of the bounded torus which has branching index $4$.

\begin{figure}[H]
\includegraphics[height=32mm, width=38mm]{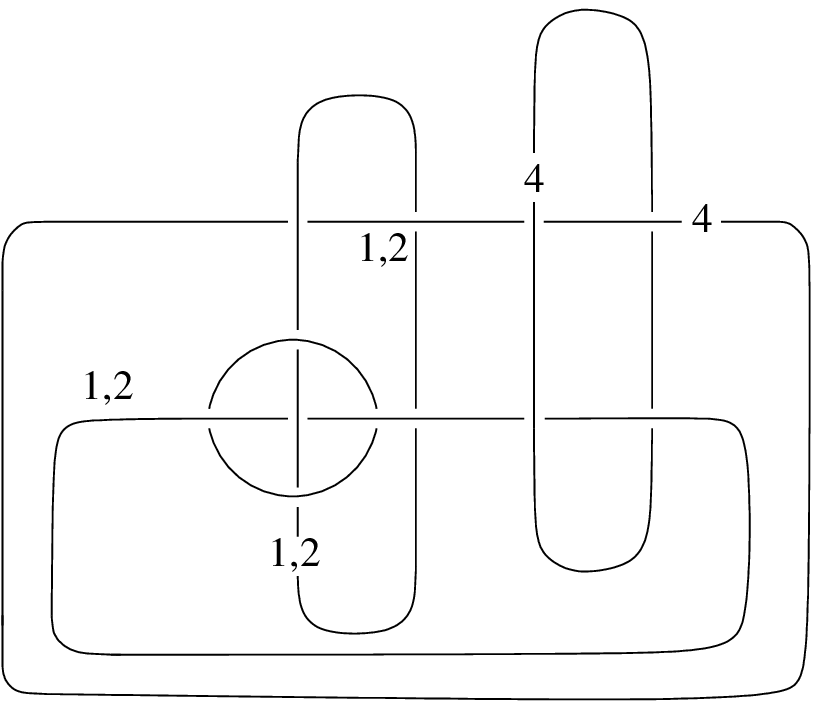}
\centering
\caption{}\label{fig:fig15}
\end{figure}

\begin{figure}[H]
\includegraphics[height=25mm, width=65mm]{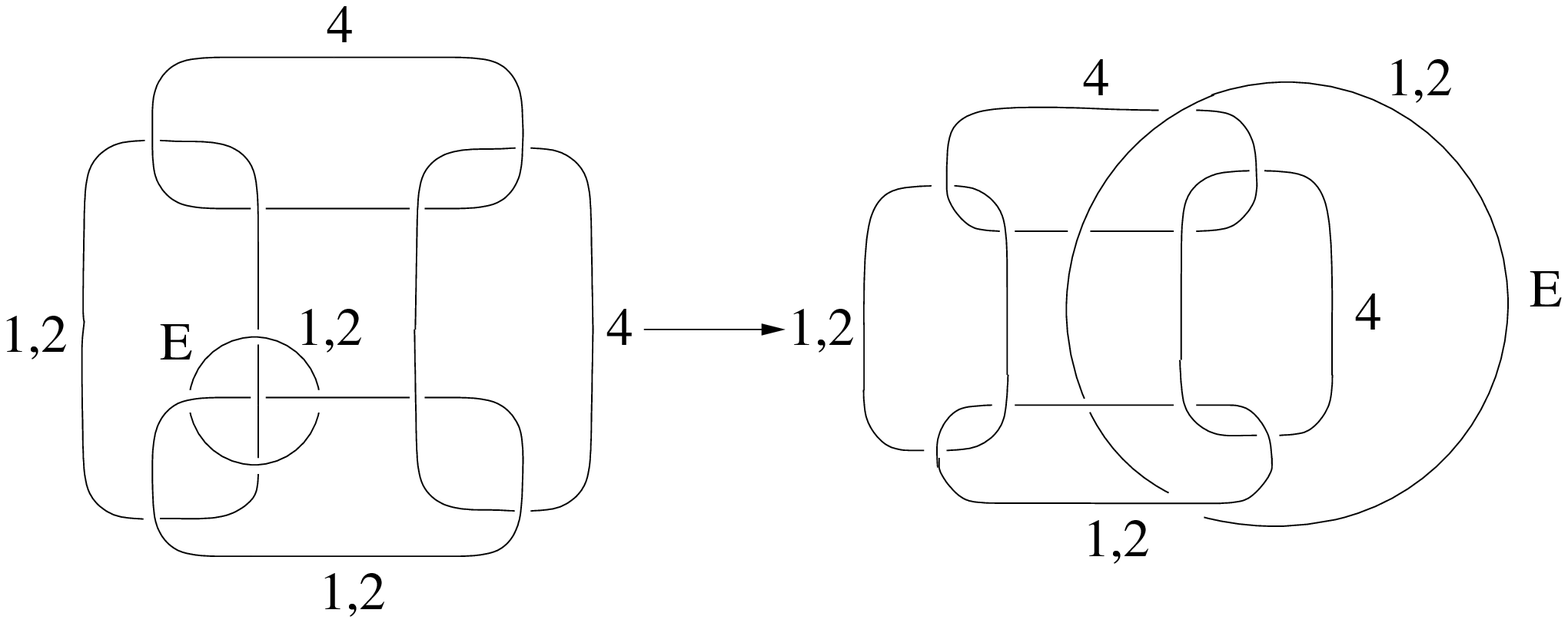}
\centering
\caption{}\label{fig:fig16}
\end{figure}

We now repeat this argument using $u$ instead of $b$. We obtain $f_4 \colon M^3 \rightarrow S^3$ with branching set $K_4$ whose branching indices are indicated in Figures \ref{fig:fig15} and \ref{fig:fig16}.

\begin{figure}[H]
\includegraphics[height=35mm, width=58mm]{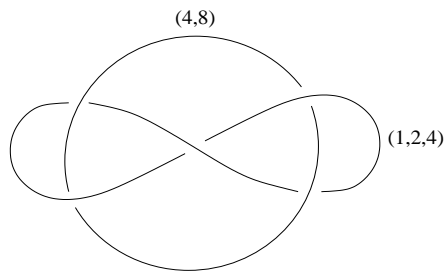}
\centering
\caption{The Whitehead link with its branching indices}\label{fig:fig17}
\end{figure}

Composing $f_4$ with the $4$-fold covering $S^3 \rightarrow S^3$ defined by a $90^\circ$-rotation about the axis $E$ of Figure \ref{fig:fig16} we obtain $f_5 \colon M^3 \rightarrow S^3$ branched over the Whitehead link with branching indices indicated in Figure \ref{fig:fig17}. This concludes the proof of Theorem 1. $\square$

We can state Theorem 1 in a different form, more suitable for future applications.

\begin{theo}
Given a closed, oriented 3-manifold M there exists a transitive representation $\omega \colon \pi_1 (S^3 - W) \rightarrow S_t$ sending the meridians $x, y$ of the Whitehead link $W$ to a product of cycles of lenghts $1,2,4$ and $4,8$ resepectively such that the manifold covering $S^3$ and branched over $W$ defined by $\omega$ is $M$. $\square$

\end{theo}

\section{The associated regular covering}

Given a transitive representation
\begin{equation}\notag
\omega \colon \pi_1 (S^3 - W) \rightarrow S_t
\end{equation}
\\
we have a $t$-fold unbranched covering of $S^3 - W$ corresponding to the stabilizer of $1$ under $\omega$, i.e., $Sb(1)=\{x\in \pi_1 (S^3 - W)\colon \omega(x)(1)=1\}$. This covering has a unique completion $f \colon M \rightarrow S^3$ which we call the branched covering defined by $\omega$ (see \cite{Fox2}).

The image of $\omega$, $G$, is a finite group of order $r$ say. It admits a canonical representation $\eta$ into $S_r$ thought of as the permutation group of the set $G$. To define $\eta$ it is sufficient to represent an element $y$ of $G$ into the bijection $\eta(y)$ of $G$ obtained by right multiplication:

\begin{align}\notag
\eta(y)\colon G \rightarrow G \\
y_i \rightarrow y_i y\notag
\end{align} 

Therefore, besides $\omega$ we have the representation $\rho = \eta \;\omega : \pi_1 (S^3 - W) \rightarrow S_r$, which gives rise to a covering $g \colon R^3 \rightarrow S^3$ branched over $W$ and corresponding to
\begin{equation}\notag
Sb(y)=\{x \in \pi_1 (S^3 - W) \colon \rho (x)(y)=y, \text{ for some fixed }y \in G\}= \text{Ker }\omega.
\end{equation}
Since this is a normal group contained in $Sb(1)$, we have a commutative diagram of branched coverings where $g$ and $h$ are regular :
\[
\xymatrix{
R^3 \ar[r]^h \ar[rd]_g &M^3 \ar[d]^f \ar[d] \\ &S^3
}
\]
Here $g$ is called the {\it associated regular covering} (it is the associated principal bundle in fibre-bundle theory) and is branched over $W$.
\\

\textbf {Addendum to Theorem 1.2.}\;{\it The branching indices of the associated regular covering defined by $\pi = \eta \;\omega$ are 4 and 8, for respective components of W.}
\\

{\it Proof.} The permutation $\eta(y), (y \in G)$ is a bijection of the set $G$, and the cycle generated by some $y_i \in G$ is 
\\
\begin{equation}\notag
(y_i ,y_iy,y_i y^2,\dotsc ,y_i y^{k-1})
\end{equation}
\\
where $k$ is the smallest positive number such that $y_i = y_i y^k$, i.e. $k$ is the order of $y$. Since $\omega(x), \omega(y)$ have orders $4$ and $8$ respectively, we conclude that $\rho(x), \rho(y)$ are products of cycles of lengths 4 and 8 respectively in $S_r$.\; $\Box$

\section {Orbifold structure of the Whitehead link}
\begin{figure}[H]
\includegraphics[height=60mm,width=80mm]{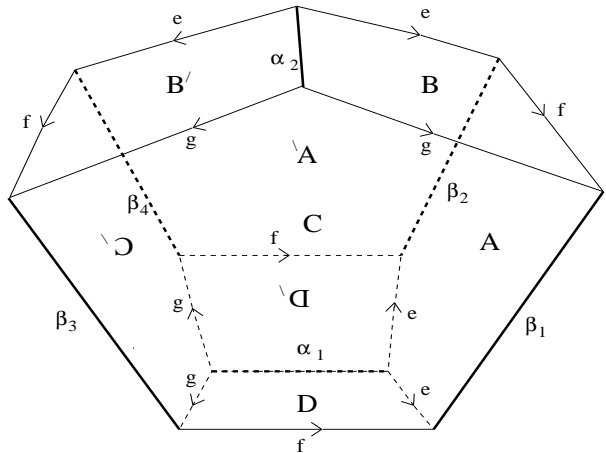}
\centering
\caption{Polyhedron with identifications for Whitehead link cone manifold}\label{fig:wlpoly}
\end{figure}

In \cite{Aal} we constructed a polyhedron with identifications as in Figure \ref{fig:wlpoly} which gives the Whitehead link inside $S^3$.

The polyhedron lives in the interior of a sphere which represents the Klein model for $\Bbb H^3$. We can tesselate $\Bbb H^3$ by this polyhedron and the group $U$ of orientation-preserving isometries of $\Bbb H^3$ generated by rotations in the respective axes defines a regular covering

\begin{equation}\notag
u:\Bbb H^3 \rightarrow S^3
\end{equation}
\\
branched over the Whitehead link $W$ with branching indices $4$ and $8$ in respective components of $W$.

\begin{figure}[H]
\includegraphics[height=30mm, width=45mm]{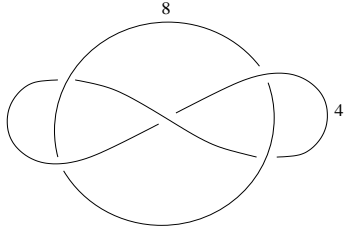}
\centering
\caption{Whitehead link cone manifold}\label{fig:fig18}
\end{figure}

{\it Remark}. The group $U$ is finitely generated. Since $U$ acts properly discontinuosly on $\Bbb H^3$, the orbit space $\Bbb H^3/U$ has the structure of an orbifold (\cite{Th2}, Proposition 13.2.1, p.13.7); in fact the orbifold of Figure \ref{fig:fig18}. Since $\Bbb H^3$ is simply connected, it follows from the remarks of Thurston (\cite{Th2}, p.13.11) that the branched cover, $u \colon \Bbb H^3 \rightarrow S^3$, is a ``universal branched cover''; that is $\Bbb H^3$ is the universal covering orbifold of $\Bbb H^3/U$. Hence, there exists an orbifold covering, $k \colon \Bbb H^3 \rightarrow R^3$, such that $gk=u$. We shall construct $k$ explicitly, in the next section.

\section{A universal group}

So far we have the coverings

\begin{align}\notag
g:R^3 \rightarrow S^3 \\ \notag
u:\Bbb H^3 \rightarrow S^3
\end{align}
\\
branched over the Whitehead link $W$ with branching indices $4$ and 8 on appropriate components. We now define a map $k:\Bbb H^3 \rightarrow R^3$ such that $gk=u$. To do this, we first orient the components $W_1$ and $W_2$ of $W$ and their preimages $g^{-1}W_i := \widetilde{W_i}$ and $u^{-1}W_i := \widehat{W_i}, i=1,2$, in such a way that both $g \colon \widetilde{W_i} \rightarrow W_i$ and $u:\widehat{W_i} \rightarrow W_i$ are orientation-preserving. Secondly, we select base points $Q, \widetilde{Q}$ and $\widehat{Q}$ in $S^3 - E^2$, $R^3$ and $\Bbb H^3$ such that $g(\widetilde{Q})=Q$ and $u(\widehat{Q})=Q$. 

Each path $\hat{\gamma}$ in $\Bbb H^3 - \widehat{W}$ starting in $\widehat{Q}$ and ending in $\widehat{Q}(\hat{\gamma})\in u^{-1}Q$ projects onto a loop $\gamma$ in $S^3 - W,$ based at $Q$. This loop lifts to a path $\tilde{\gamma}$ starting in $\widetilde{Q}$ and ending in $\widetilde{Q}(\tilde{\gamma})\in g^{-1}Q.$
We map the polyhedron $P(\hat{\gamma})$ containing $\widehat{Q}(\hat{\gamma})$ onto the closure $P(\tilde{\gamma})$ of the component of $g^{-1}(S^3 - E^2)$ containing $\widetilde{Q}(\tilde{\gamma})$ recalling that $E^2$ is a concrete splitting complex for $W \subset S^3$. The map

\begin{equation}\notag
k(\tilde{\gamma})\colon P(\hat{\gamma}) \rightarrow P(\tilde{\gamma})
\end{equation}
\\
is a simplicial map: it sends vertices, edges and faces of $P(\hat{\gamma})$ into vertices, edges and faces of $P(\tilde{\gamma})$ in such a way that $\widehat{W_i} \cap P(\hat{\gamma})$ maps onto $\widetilde{W_i} \cap P(\tilde{\gamma})$ orientation preservingly for $i=1,2$.

The map $k(\tilde{\gamma})$ is independent of $\hat{\gamma}$. For if $\hat{\gamma}' \subset \Bbb H^3 - \widehat{W}$ also joins $\widehat{Q}$ and $\widehat{Q} (\hat{\gamma})$, then $\hat{\gamma} * (\hat{\gamma}')^{-1}$ is homotopic in $\Bbb H^3 - \widehat{W}$ to a product of meridians $m_1 * m_2 * \dots * m_s$ of $\widehat{W}$. Therefore, $\gamma * (\gamma')^{-1}$ is homotopic in $S^3 - W$ to a product $M_1^{n_1} * \dots * M_s^{n_s}$ where $M_1,\dots,M_s$ are meridians of $W$ and $n_1,\dots ,n_s \in \{4,8\}$. Thus, since $\rho(M_i)$ has order 4 or 8, $\rho(\gamma)=\rho(\gamma')$. Hence $\tilde{\gamma}$ and $\tilde{\gamma}'$ end in the same point of $R^3 - \widetilde{W}$.

\begin{theo}The map $k \colon \Bbb H^3 \rightarrow R^3$ is the universal covering of $R^3$ and \\$gk=u$.
\end{theo}

{\it Proof} \; By construction $gk=u$, and $k | \Bbb H^3 - \widehat{W}$ is a covering from $\Bbb H^3 - \widehat{W}$ onto $R^3 - \widehat{W}$. Since the branching indices of $g$ and $u$ are 4 or 8, in both cases, it follows that $k$ is a covering. $\Box$ 

\begin{cor} The group U, generated by $\pi/2$ and $\pi/4$ rotations in the appropriate axes of the polyhedron $P$ is universal, i.e. every closed, oriented 3-manifold $M^3$ is orientation-preservingly homeomorphic to $\Bbb H^3/G$, where $G$ is a subgroup of $U$ of finite index. Thus, every manifold $M^3$ is a hyperbolic orbifold covering the hyperbolic orbifold with underlying space $S^3$ and singular locus the Whitehead link with isotropy groups $\Bbb Z/4\Bbb Z$ and $\Bbb Z/8\Bbb Z$ on respective link components.
\end{cor} 

{\it Proof} \; In general, given a branched covering $p \colon M \rightarrow N$ we have the group Aut$(M,p)$ of fibre-preserving automorphisms of $M$. In our case $U$ equals Aut$(\Bbb H^3,u)$ and $\Bbb H^3/U$ is $S^3$ so that the branched covering $u \colon \Bbb H^3 \rightarrow S^3$ is equivalent to the natural projection $\Bbb H^3 \rightarrow \Bbb H^3/U$. Thus $u$ is a regular covering. Therefore, since $u=f\cdot(hk)$, it follows that $hk \colon \Bbb H^3 \rightarrow M^3$ is regular and that Aut$(\Bbb H^3,hk)$ is a subgroup $G$ of $U$. Moreover, the index of $G$ is the number of sheets of $f \colon M \rightarrow S^3$. $\Box$

\begin{cor} If the closed oriented 3-manifold $M$ is homeomorphic to $\Bbb H^3/G$, where $G$ is a finite index subgroup of $U$, then $\pi_1(M)$ is isomorphic to $G/F$, where $F$ is the normal subgroup of $G$ generated by the elements with fixed points (i.e. the elements of finite order).
\end{cor}

{\it Proof} \; Consequence of \cite{Arm}. $\Box$

\begin{cor} $M \cong \Bbb H^3/G$ is simply connected if and only if $G$ is generated by elements of finite order. $\Box$
\end{cor}

Later we will show that $W(4n,8n)$ is hyperbolic for $0\leq n \leq 1$ and we will study the topology of $\Bbb H^3/G$ for finitely generated finite index $G\leq U=W(4,8)$.

\bibliographystyle{amsplain}
\bibliography{pc1}

\vspace{6mm}
\hspace{-5mm}Dr Aalam\\
PO Box 18810\\
London SW7 2ZR\\
email : aalam@mth.kcl.ac.uk\\ali\_aalam@hotmail.com

\end{document}